%% file: lhpp_groupoids.tex
\documentclass[12pt,a4paper]{article}
\pdfoutput=1
\usepackage{amsmath}
\usepackage{amsfonts}
\usepackage{amssymb}
\usepackage{amsthm}
\usepackage{graphicx}
\usepackage[all]{xy}

\newtheorem{thm}{Theorem}[section]
\newtheorem{prop}{Proposition}[section]
\newtheorem{dfn}{Definition}[section]

\theoremstyle{remark}
\newtheorem{rmk}{Remark}[section]
\newtheorem{expl}{Example}[section]

\newcommand{\R}{\mathbb{R}}
\newcommand{\categories}[1]{\mathbf{#1}}
\newcommand{\liebr}[2]{[{#1},{#2}]}

\newcommand{\T}{\mathbf{T}}
\newcommand{\Pb}[4]{#1\;{}_{#2}\!\times_{#3}#4}
\DeclareMathOperator{\Aut}{Aut}
\DeclareMathOperator{\End}{End}
\DeclareMathOperator{\Iso}{Iso}
\DeclareMathOperator{\Hom}{Hom}

\author{Luiz Henrique P. P\^egas\footnote{Contact by e-mail: lhp@ime.usp.br or divtzero@gmail.com}}
\title{Groupoids and cogroupoids: an one object approach}
\date{\today}

\begin{document}
\maketitle

\begin{abstract}
The aim of this paper is to provide a definition of groupoid and cogroupoid internal to a category which makes use of only one object and morphisms, in contrast with the two object approach commonly found in the literature. We will give some examples and we will stablish a relation with group objects (and Hopf algebras). The definitions presented here were designed to simplify some constructions related to internalizations, Lie groupoids and Hopf algebroids.
\end{abstract}

\tableofcontents

\input{body_text}

\appendix
\input{appendix}

\bibliographystyle{plain}
\bibliography{lhpp_groupoids}
\end{document}

%% file: body_text.tex
\section{Introduction}
In standard texts about groupoids is often emphasized the bundle nature of groupoids, even those in which the approach is mainly categorial. However, a groupoid is an algebraic object in its very essence and such a point of view is capable to clarify some simpler aspects that may be hidden by the complexity of the bundle structure. In these notes we present an algebraic point of view which has some advantages in some situations (and perhaps some disadvantages in others).

In what follows we intend to use diagrams to express properties, meaning, axioms fulfilled in a given context. In this spirit, when an object is mentioned, the reader can mentally add an ``if exists'' in the argument.

About the notation, given a category $\categories{C}$, the class of objects will be denoted by $\categories{C}_{0}$, the class of morphisms by $\categories{C}_{1}$, the class of morphisms from the object $A$ to the object $B$ by $\categories{C}(A,B)$ and $\categories{C}(A,A)$ by $\End_{\categories{C}}(A)$. The tangent functor on smooth manifolds will be denoted by $\T$. It may not be usual, but the choice was made to avoid conflict with the uppercase Greek letter tau ($T$), since we need the symbols $\Sigma$ and $T$ to represent the source and target projections, in a slightly different meaning of the standard $\sigma$ and $\tau$. Sometimes, we will use the circle notation for composites, but often it will be denoted just by juxtaposition.

We will deal with pull-backs of morphisms which are often ``legs'' of spans. When we want this emphasis, for the pull-back of $A' \leftarrow A \rightarrow C$ and $C~\leftarrow~B~\rightarrow~B'$, instead of writing
\begin{displaymath}
\xymatrix{\Pb{A}{f}{g}{B} \ar[d]_-{\pi_{fg1}} \ar[r]^-{\pi_{fg2}} & B \ar[d]^-{g} \\
A \ar[r]_-{f} & C}
\end{displaymath}
we will write
\begin{displaymath}
\xymatrix{ & & \Pb{A}{f}{g}{B} \ar[dl]_-{\pi_{fg1}} \ar[dr]^-{\pi_{fg2}} & & \\
& A \ar[dl]_-{f'} \ar[dr]_-{f} & & B \ar[dl]^-{g} \ar[dr]^-{g'} & \\
A' & & C & & B'}
\end{displaymath}
obtaining a new span $A' \xleftarrow{f' \circ \pi_{fg1}} \Pb{A}{f}{g}{B} \xrightarrow{g' \circ \pi_{fg2}} B'$.

In section \ref{groupoid_struc_sec} we will discuss the structure of groupoids in a category, with examples of groupoids in smooth manifolds (Lie groupoids) and groupoids in associative algebras. The section ends with the obtention of group objects as a special case of groupoid objects.

In section \ref{cogroupoids_sec} we will discuss the structure of cogroupoids in a category, with an example in commutative associative unital algebras and the obtention of Hopf algebras as a special case of cogroupoids in that algebras.

In section \ref{groupoid_actions_sec} we will discuss a diagrammatic description of an action of a groupoid. The section ends with an application to groupoids in smooth manifolds: the construction of the Lie algebroid of a Lie groupoid by means of projections. The last section only indicates directions of work by this author. In the appendix there are some remarks that do not fit well in the body of the text.

\section{The groupoid structure}
\label{groupoid_struc_sec}
The classical definition of a groupoid (see \emph{e.g.} \cite{mackenzie_lie_groupoids}) makes extensive use of a (double) bundle structure. It is possible to avoid a direct use of that structure, encoding all that information in a single object sketch. Let us start with the sketch of a category in a category.
\begin{dfn}[Category in a category]
\label{category_dfn}
Let $\categories{C}$ be a category, $C \in \categories{C}_{0}$, $\Sigma,T \in \End_{\categories{C}}(C)$ and $\mu \in \categories{C}(C_{2},C)$, where $C_{2} \in \categories{C}_{0}$ is the pull-back
\begin{equation}
\label{Cat2_def}
\xymatrix{ & & \Pb{C}{\Sigma}{T}{C} \ar[dl]_-{\pi_{\Sigma T1}} \ar[dr]^-{\pi_{\Sigma T2}} & & \\
& C \ar[dl]_-{T} \ar[dr]_-{\Sigma} & & C \ar[dl]^-{T} \ar[dr]^-{\Sigma} & \\
C & & C & & C}
\end{equation}
We call $(C,\Sigma,T,\mu)$ a category in $\categories{C}$ if these data satisfy
\begin{enumerate}
\item[i)] $T \Sigma = \Sigma$ and $\Sigma T=T$;
\item[ii)] it is commutative
\begin{equation}
\label{mu_tau_sigma_def}
\xymatrix{ & \Pb{C}{\Sigma}{T}{C} \ar[dl]_-{T \circ \pi_{\Sigma T1}} \ar[dd]^-{\mu} \ar[dr]^-{\Sigma \circ \pi_{\Sigma T2}} & \\
C & & C \\
& C \ar[ul]^-{T} \ar[ur]_-{\Sigma} & }
\end{equation}
\item[iii)] by denoting $\eta_{\Sigma} \colon C \to \Pb{C}{\Sigma}{T}{C}$ as the unique arrow given by the cone $C \xleftarrow{id_{C}} C \xrightarrow{\Sigma} C$ such that $\pi_{\Sigma T1} \circ \eta_{\Sigma}=id_{C}$ and $\pi_{\Sigma T2} \circ \eta_{C}=\Sigma$, and similarly $\eta_{T} \colon C \to \Pb{C}{\Sigma}{T}{C}$ the arrow given by the cone $C \xleftarrow{T} C \xrightarrow{id_{C}} C$ such that $\pi_{\Sigma T1} \circ \eta_{T}=T$ and $\pi_{\Sigma T2} \circ \eta_{T}=id_{C}$, then
\begin{eqnarray}
\mu \circ \eta_{\Sigma} & = & id_{C} \nonumber \\
\label{eta_S_mu_def}\mu \circ \eta_{T} & = & id_{C}
\end{eqnarray}
\item[iv)] there exist the pull-back $\Pb{(\Pb{C}{\Sigma}{T}{C})}{\Sigma_{2}}{T}{C}$, given by
\begin{displaymath}
\xymatrix{
& & \Pb{(\Pb{C}{\Sigma}{T}{C})}{\Sigma_{2}}{T}{C} \ar[dl]_-{\pi_{\Sigma_{2}T1}} \ar[dr]^-{\pi_{\Sigma_{2}T2}} & & \\
& \Pb{C}{\Sigma}{T}{C} \ar[dl]_-{T \circ \pi_{\Sigma T1}} \ar[dr]^-{\Sigma \circ \pi_{\Sigma T2}} & & C \ar[dl]_-{T} \ar[dr]^-{\Sigma} & \\
C & & C & & C}
\end{displaymath}
and the pull-back $\Pb{C}{\Sigma}{T_{2}}{(\Pb{C}{\Sigma}{T}{C})}$ given by a similar diagram;
\item[v)] by denoting $\mu \times id_{C} \colon \Pb{(\Pb{C}{\Sigma}{T}{C})}{\Sigma_{2}}{T}{C} \to \Pb{C}{\Sigma}{T}{C}$ as the unique arrow given by the diagram
\begin{displaymath}
\xymatrix{
\Pb{(\Pb{C}{\Sigma}{T}{C})}{\Sigma_{2}}{T}{C} \ar@/^1pc/[drr]^-{id_{C} \circ \pi_{\Sigma_{2}T2}} \ar@{-->}[dr]^-{\mu \times id_{C}} \ar@/_1pc/[ddr]_-{\mu \circ \pi_{\Sigma_{2}T1}} & & \\
& \Pb{C}{\Sigma}{T}{C} \ar[r]_-{\pi_{\Sigma T2}} \ar[d]^-{\pi_{\Sigma T1}} & C \ar[d]^-{T} \\
& C \ar[r]_-{\Sigma} & C}
\end{displaymath}
which is well defined, because
\begin{displaymath}
\Sigma \circ \mu \circ \pi_{\Sigma_{2}T1} = \Sigma \circ \pi_{\Sigma T2} \circ \pi_{\Sigma_{2}T1} = T \circ \pi_{\Sigma_{2}T2} = T \circ id_{C} \circ \pi_{\Sigma_{2}T2}
\end{displaymath}
and by $id_{C} \times \mu \colon \Pb{C}{\Sigma}{T_{2}}{(\Pb{C}{\Sigma}{T}{C})} \to \Pb{C}{\Sigma}{T}{C}$ the unique arrow given by the analogous cone, it holds true
\begin{equation}
\label{mu_associative_def}
\mu \circ (\mu \times id_{C}) \approx \mu \circ (id_{C} \times \mu)
\end{equation}
meaning that the left hand side of the last equation is the same as the right hand side up to an isomorphism. This is necessary because we have in general an isomorphism between $\Pb{(\Pb{C}{\Sigma}{T}{C})}{\Sigma_{2}}{T}{C}$ and $\Pb{C}{\Sigma}{T_{2}}{(\Pb{C}{\Sigma}{T}{C})}$ and not the equality;
\item[vi)] the pairs $(\Sigma,id_{C})$ and $(T,id_{C})$ have coequalisers.
\end{enumerate}
\end{dfn}
\begin{rmk}
There is an important consequence of the conditions in $\emph{i)}$ above, with a strong algebraic flavour. From this two conditions, we can see that
\begin{displaymath}
\Sigma \Sigma=\Sigma
\end{displaymath}
\begin{displaymath}
T T=T
\end{displaymath}
These are simple calculations: for example,
\begin{displaymath}
T T=T \Sigma T=\Sigma T=T
\end{displaymath}
The other equation is similar. Hence, conditions $\emph{i)}$ and $\emph{vi)}$ are statements that $\Sigma$ and $T$ are split idempotents. $\Diamond$
\end{rmk}
\begin{dfn}[Groupoid in a category]
\label{groupoid_dfn}
Let $\categories{C}$ be a category, $(G,\Sigma,T,\mu)$ a category in $\categories{C}$ and $\Upsilon \in \End_{\categories{C}}(G)$. We call $(G,\Sigma,T,\Upsilon,\mu)$ a groupoid in $\categories{C}$ if these satisfy
\begin{enumerate}
\item[i)] $T \Upsilon=\Sigma$, $\Upsilon \Sigma=\Sigma$, $\Upsilon \Upsilon=id_{G}$;
\item[ii)] by denoting $id_{G} \times \Upsilon \colon G \to \Pb{G}{\Sigma}{T}{G}$ as the unique arrow given by the diagram
\begin{displaymath}
\xymatrix{
G \ar@/^1pc/[drr]^-{\Upsilon} \ar@{-->}[dr]^-{id_{G} \times \Upsilon} \ar@/_1pc/[ddr]_-{id_{G}} & & \\
& \Pb{G}{\Sigma}{T}{G} \ar[r]_-{\pi_{\Sigma T2}} \ar[d]^-{\pi_{\Sigma T1}} & G \ar[d]^-{T} \\
& G \ar[r]_-{\Sigma} & G}
\end{displaymath}
and by $\Upsilon \times id_{G} \colon G \to \Pb{G}{\Sigma}{T}{G}$ the analogous arrow, which is well defined, because
\begin{displaymath}
\Sigma \Upsilon = T \Sigma \Upsilon = T T \Upsilon \Upsilon = T \circ id_{G}
\end{displaymath}
it holds true both
\begin{eqnarray}
\mu \circ (id_{G} \times \Upsilon) & = & T \nonumber \\
\label{mu_upsilon_eqn}\mu \circ (\Upsilon \times id_{G}) & = & \Sigma
\end{eqnarray}
\end{enumerate}
\end{dfn}
\begin{rmk}
The definition of groupoid given here can be rephrased as: a groupoid in $\categories{C}$ is an internal category in $\categories{C}$ such that every morphism is an isomorphism, where category is defined as an internalisation of the definition of category given in \cite{freyd_categories}. So, in essence, there is nothing new about this definition. The usual notion of groupoid can be easily recovered by noting that there is a bijection between the class of identities in a category and the class of objects. In the internal sense, this is achieved by means of the coequalisers. $\Diamond$
\end{rmk}
The equivalence between this definition of internal category and the usual one can be codified in the following proposition.
\begin{prop}
\label{object_of_objects_prop}
Let $\categories{C}$ be a category with the necessary products and $(C,\Sigma,T,\mu)$ a category in $\categories{C}$. Then the statement that the pairs $(\Sigma,id_{C})$ and $(T,id_{C})$ have coequalisers is the same that stating there is a subobject $M$ of $C$ such that the diagram
\begin{displaymath}
\xymatrix@C+50pt{M \ar[r]^-{\varepsilon} \ar[d]_-{\varepsilon} & C \ar[d]^-{(id_{C},\Sigma,T)} \\
C \ar[r]_-{(id_{C},id_{C},\Sigma)} & C \times C \times C}
\end{displaymath}
is a pull-back diagram.
\begin{proof}
Suppose that $\sigma \colon C \to M$ and $\tilde{\tau} \colon C \to N$ are coequalizers for $(\Sigma,id_{C})$ and $(T,id_{C})$ respectively. By the fact that $\Sigma \Sigma=\Sigma$, there exists a unique $\varepsilon \colon M \to C$ such that $\varepsilon \sigma=\Sigma$. Using again the universal property of the coequalizer, is straightforward that $\sigma \varepsilon=id_{M}$. We can apply the same reasoning to $\tilde{\tau}$, giving a unique $\tilde{\varepsilon} \colon N \to C$, such that $\tilde{\varepsilon} \tilde{\tau}=T$ and $\tilde{\tau} \tilde{\varepsilon}=id_{N}$. Now we claim that $M$ and $N$ are isomorphic, with isomorphisms $\sigma \tilde{\varepsilon}$ and $\tilde{\tau} \varepsilon$, mutually inverse. To see this, just notice that
\begin{eqnarray*}
& & \Sigma T = T \\
& & \varepsilon \sigma \tilde{\varepsilon} \tilde{\tau} = \tilde{\varepsilon} \tilde{\tau} \\
& & \varepsilon \sigma \tilde{\varepsilon} = \tilde{\varepsilon} \\
& & \tilde{\tau} \varepsilon \sigma \tilde{\varepsilon} = id_{N}
\end{eqnarray*}
and analogously, from $T \Sigma=\Sigma$, that $\sigma \tilde{\varepsilon} \tilde{\tau} \varepsilon = id_{M}$. It is worth noticing that $\varepsilon$ and $\tilde{\varepsilon}$ are monomorphisms, because if we have $f$ and $g$ such that $\varepsilon f = \varepsilon g$, then $\sigma \varepsilon f = \sigma \varepsilon g$ implies $f=g$ (and the same for $\tilde{\varepsilon}$).

By defining $\tau=\sigma \tilde{\varepsilon} \tilde{\tau}$, we have $\tau \colon C \to M$, such that
\begin{eqnarray*}
& & \varepsilon \tau = \varepsilon \sigma \tilde{\varepsilon} \tilde{\tau} = \Sigma T = T \\
& & \tau \varepsilon = \sigma \tilde{\varepsilon} \tilde{\tau} \varepsilon = id_{M}
\end{eqnarray*}
Hence, we have
\begin{displaymath}
\Sigma \varepsilon = \varepsilon \sigma \varepsilon = \varepsilon
\end{displaymath}
\begin{displaymath}
\Sigma \varepsilon = \varepsilon = \tilde{\varepsilon} \tilde{\tau} \varepsilon = T \varepsilon
\end{displaymath}
where we used that $\varepsilon = \tilde{\varepsilon} \tilde{\tau} \varepsilon$ (which is a consequence of $T \Sigma = \Sigma$). Hence, we proved that the diagram
\begin{displaymath}
\xymatrix@C+50pt{M \ar[r]^-{\varepsilon} \ar[d]_-{\varepsilon} & C \ar[d]^-{(id_{C},\Sigma,T)} \\
C \ar[r]_-{(id_{C},id_{C},\Sigma)} & C \times C \times C}
\end{displaymath}
is commutative. To show that it is a pull-back, let $C \xleftarrow{\beta} D \xrightarrow{\alpha} C$ be a cone such that $(id_{C},\Sigma,T) \circ \alpha = (id_{C},id_{C},\Sigma) \circ \beta$. Notice that we have $\alpha=\beta$. Define $\eta \colon D \to M$ by $\eta=\sigma \alpha$. Thus,
\begin{displaymath}
\varepsilon \eta = \varepsilon \sigma \alpha = \Sigma \alpha = \alpha
\end{displaymath}
Now, $\eta$ is the unique arrow with this property, because $\varepsilon$ is monomorphism.

By the other hand, suppose that the given diagram is a pull-back. First of all, note that $\varepsilon \colon M \to C$ is an equalizer for $(\Sigma,id_{C})$. To show this, note that $\Sigma \varepsilon=\varepsilon$ and if $e \colon N \to C$ is another arrow such that $\Sigma e=e$, then
\begin{displaymath}
T e = T \Sigma e = \Sigma e = e
\end{displaymath}
By the universal property of the pull-back, there exists a unique $\bar{e} \colon N \to M$ such that $\varepsilon \bar{e}=e$, which gives the universal property of the equalizer. Thus, by the fact that $\Sigma \Sigma=\Sigma$, there exists a unique arrow $\sigma \colon C \to M$ such that $\varepsilon \sigma = \Sigma$. Besides, such an arrow is a left inverse for $\varepsilon$. To see this, notice that
\begin{displaymath}
\varepsilon \sigma \varepsilon = \Sigma \varepsilon = \varepsilon
\end{displaymath}
By the universal property of the equalizer, any arrow $f$ such that $\varepsilon f = \varepsilon$ is $f=id_{M}$. Therefore, $\sigma \varepsilon = id_{M}$. Now, $\sigma$ is a coequalizer for $(\Sigma,id_{C})$. Indeed,
\begin{displaymath}
\sigma \Sigma = \sigma \varepsilon \sigma = \sigma
\end{displaymath}
and if $s \colon C \to N$ is such that $s \Sigma = s$, then the arrow $s \varepsilon \colon M \to N$ satisfies
\begin{displaymath}
s \varepsilon \sigma = s \Sigma = s
\end{displaymath}
Suppose that $\bar{s} \colon M \to N$ is an arrow such that $\bar{s} \sigma = s$. Then
\begin{displaymath}
\bar{s} = \bar{s} \sigma \varepsilon = s \varepsilon = s \varepsilon \sigma \varepsilon = s \varepsilon
\end{displaymath}
showing the result. A similar argument shows that there is a unique $\tau \colon C \to M$, with $\tau \varepsilon = id_{M}$ and $\varepsilon \tau = T$, such that $\tau$ is a coequalizer for $(T,id_{C})$.
\end{proof}
\end{prop}
\begin{rmk}
\label{product_delta_rmk}
There are some cases (some categories) in which pull-backs can be described as products followed by equalizers. In these situations we have a ``factorization'' of some arrows, which can be useful. Let $(G,\Sigma,T,\Upsilon,\mu)$ be a groupoid in a category with such factorization and denote by $\Delta \colon G \to G \times G$ the unique arrow given by the cone $G \xleftarrow{id_{G}} G \times G \xrightarrow{id_{G}} G$. Then, we can write
\begin{eqnarray*}
\eta_{\Sigma} & = & (id_{G},\Sigma) \circ \Delta \\
\eta_{T} & = & (T,id_{G}) \circ \Delta \\
\Upsilon \times id_{G} & = & (\Upsilon,id_{G}) \circ \Delta \\
id_{G} \times \Upsilon & = & (id_{G},\Upsilon) \circ \Delta
\end{eqnarray*} $\Diamond$
\end{rmk}
The first advantage of definition \ref{groupoid_dfn} can be seen when we are dealing with homomorphisms of groupoids.
\begin{dfn}[Homomorphisms of groupoids]
\label{homomorphism_groupoids_dfn}
Let $\categories{C}$ be a category and $(G,\Sigma_{G},T_{G},\Upsilon_{G},\mu_{G})$, $(H,\Sigma_{H},T_{H},\Upsilon_{H},\mu_{H})$ be groupoids in $\categories{C}$. $F \in \categories{C}(G,H)$ is called a homomorphism of groupoids (or a homomorphism of the groupoids $G$ and $H$), if and only if the following diagrams are commutative
\begin{displaymath}
\xymatrix{G \ar[d]_{\Sigma_{G}} \ar[r]^{F} & H \ar[d]^{\Sigma_{H}} & & G \ar[d]_{T_{G}} \ar[r]^{F} & H \ar[d]^{T_{H}} & & G_{2} \ar[d]_{\mu_{G}} \ar[r]^{F \times F} & H_{2} \ar[d]^{\mu_{H}} \\
G \ar[r]_{F} & H & & G \ar[r]_{F} & H & & G \ar[r]_{F} & H}
\end{displaymath}
\end{dfn}
\begin{rmk}
Notice that if the above diagrams are commutative, so is the corresponding diagram involving $\Upsilon_{G}$ and $\Upsilon_{H}$. It is also worth noticing that this is just the definition of internal functor between categories ``only with morphisms''. $\Diamond$
\end{rmk}
With these definitions at hand, it is possible to discuss the concept of groupoid in any category with some products and some pull-backs. As one can easily see, when in a concrete category, the double bundle structure is recovered by the coequalizers conditions, which give the base object. The reader is invited to do some sketches in categories like $\categories{Set}$, the category of functions between (ZFC) sets, or $\categories{Top}$, the category of continuous functions between (small) topological spaces.

The above definition was designed to shed some light on facts about groupoids in $\categories{Man}$, the category of smooth ($C^{\infty}$) functions between smooth manifolds (locally Euclidean, second countable, Hausdorff topological spaces).\footnote{For those concerned with size problems, see appendix \ref{appendix_sec}.} The next proposition will play a central role on the discussion about Lie groupoids.
\begin{prop}
\label{embedded_base_prop}
Suppose that $(G,\Sigma,T,\Upsilon,\mu)$ is a groupoid in $\categories{Man}$. Then there is an embedded submanifold $M$ of $G$ such that $M=\Sigma(G)=T(G)$.
\begin{proof}
Suppose first $G$ connected and consider the set $M=im(\Sigma)$. If $x \in M$, by $\Sigma \Sigma = \Sigma$ we have $\T_{x}\Sigma \circ \T_{x}\Sigma = \T_{x}\Sigma$ which leads to $im(\T_{x}\Sigma)=ker(Id-\T_{x}\Sigma)$. By the kernel and image theorem, we have
\begin{displaymath}
dim(im(\T_{x}\Sigma))+dim(ker(\T_{x}\Sigma))=dim(G)
\end{displaymath}
which means
\begin{displaymath}
dim(ker(Id-\T_{x}\Sigma))+dim(ker(\T_{x}\Sigma))=dim(G)
\end{displaymath}
By continuity, there is an open neighbourhood $U_{x}$ of $x$ such that neither $dim(ker(\T_{x}\Sigma))$ nor $dim(ker(Id-\T_{x}\Sigma))$ can falls down. Hence, in that neighbourhood, $\T_{x}\Sigma$ has constant rank. From $G$ connected and $\Sigma$ continuous, it follows that $rank(\T_{x}\Sigma)=r$ for all $x \in im(\Sigma)$. For $x \in im(\Sigma)$ we have yet an open neighbourhood $V_{x}$ of $x$ such that $rank(\T_{y}\Sigma) \geq r,\;\forall y \in V_{x}$. On the other hand,
\begin{displaymath}
rank(\T_{y}\Sigma)=rank(\T_{y}(\Sigma \Sigma))=rank(\T_{\Sigma(y)}\Sigma \circ \T_{y}\Sigma) \leq rank(\T_{\Sigma(y)}\Sigma)=r
\end{displaymath}
therefore $rank(\T_{y}\Sigma)=r,\;\forall y \in V_{x}$.

It follows that $V=\bigcup_{x \in im(\Sigma)}V_{x}$ is a neighbourhood of $M=im(\Sigma)$ such that $\T\Sigma$ has constant rank. By the theorem of the constant rank map, $M$ is an embedded submanifold of $G$.

If $G$ is not connected, the argument can be repeated in each component and the result still holds.

Now, it remains to show that $im(\Sigma)=im(T)$. On one hand, we have
\begin{multline*}
x \in im(\Sigma) \Rightarrow (\exists g \in G;\; \Sigma(g)=x) \Rightarrow T\Sigma(g)=T(x) \Rightarrow \Sigma(g)=T(x) \Rightarrow \\
\Rightarrow x=T(x) \Rightarrow im(\Sigma) \subseteq im(T)
\end{multline*}
because $T\Sigma=\Sigma$. On the other hand,
\begin{multline*}
y \in im(T) \Rightarrow (\exists g \in G;\; T(g)=y) \Rightarrow \Sigma T(g)=\Sigma(y) \Rightarrow T(g)=\Sigma(y) \Rightarrow \\
\Rightarrow y=\Sigma(y) \Rightarrow im(T) \subseteq im(\Sigma)
\end{multline*}
because $\Sigma T=T$. Hence, $im(T)=im(\Sigma)=M$.
\end{proof}
\end{prop}
\begin{rmk}
The embedding $\varepsilon \colon M \to G$ is proper. To see this, let $K$ be a compact set in $G$. If $x \in \varepsilon^{-1}(K)$, then $x=T(\varepsilon(x))$ gives $x \in T(K)$, hence $\varepsilon^{-1}(K) \subseteq T(K)$. By continuity of $T$, $T(K)$ is a compact subset of $M$. By continuity of $\varepsilon$ and $K$ being closed, it follows that $\varepsilon^{-1}(K)$ is a closed subset of a compact set (in a Hausdorff space), therefore compact. This shows that $M$ is a properly (closed) embedded submanifold of $G$. $\Diamond$
\end{rmk}

The above remark is not so surprising since, in essence, is a restatement of the fact that if $(G,\Sigma,T,\Upsilon,\mu)$ is a groupoid in $\categories{Man}$, then there is an object $M$ defined by the pull-back
\begin{displaymath}
\xymatrix@C+50pt{M \ar[r]^-{\varepsilon} \ar[d]_-{\varepsilon} & G \ar[d]^-{(id_{G},\Sigma,T)} \\
G \ar[r]_-{(id_{G},id_{G},\Sigma)} & G \times G \times G}
\end{displaymath}
For the $\categories{Man}$ category, no natural description of subobjects is available, so the existence of the proper embedded manifold $M$ cannot be derived by pure categorial reasoning. The interesting thing is that the condition on coequalizers, in this case, is automatically fulfilled. The subobject $M$ in the category $\categories{Man}$ (as in many others) is uniquely determined (up to isomorphism) by the algebraic conditions on $\Sigma$ and $T$, being redundant the requirement of a ``unity section''.

The usual definition of Lie groupoids\footnote{See \emph{e.g.} \cite{mackenzie_lie_groupoids}.} involves a pair of manifolds\footnote{I always consider the base space and the total space of a Lie groupoid as being Hausdorff, second countable manifolds, but in other texts this may not be the case. For those important cases where the groupoid is able to be non Hausdorff, I prefer the terminology differentiable groupoid instead of Lie groupoid.} $G$ and $M$, a pair of surjective submersions $\sigma,\tau \colon G \to M$ (the source and the target projections), a smooth groupoid multiplication $\mu \colon G_{2} \to G$, a smooth inversion map $\Upsilon \colon G \to G$ and an embedding $\varepsilon \colon M \to G$ (the identity section), satisfying some desired conditions. The next theorem states that such a Lie groupoid is just a groupoid in $\categories{Man}$.
\begin{thm}
\label{Lie_Man_groupoid_equiv_thm}
Let $(G,\Sigma,T,\Upsilon,\mu)$ be a groupoid in $\categories{Man}$. Then there exist an embedding $\varepsilon \colon M \to G$, a pair of surjective submersions $\sigma,\tau \colon G \to M$, a smooth map $\mu \colon G_{2} \to G$ and a smooth map $\imath \colon G \to G$ such that
\begin{enumerate}
\item[i)] for $(h,g) \in G_{2}$, $\sigma(\mu(h,g))=\sigma(g)$ and $\tau(\mu(h,g))=\tau(h)$;
\item[ii)] $\mu$ is associative;
\item[iii)] $\sigma(\varepsilon(x))=\tau(\varepsilon(x))=x$ for all $x \in M$;
\item[iv)] $\mu(g,\varepsilon(\sigma(g)))=g$ and $\mu(\varepsilon(\tau(g)),g)=g$, for all $g \in G$;
\item[v)] for $g \in G$, $\sigma(\imath(g))=\tau(g)$, $\tau(\imath(g))=\sigma(g)$ and $\mu(\imath(g),g)=\varepsilon(\sigma(g))$, $\mu(g,\imath(g))=\varepsilon(\tau(g))$.
\end{enumerate}
\begin{proof}
From proposition \ref{embedded_base_prop}, the set $M=im(\Sigma)=im(T)$ is an embedded manifold in $G$. Let $\varepsilon \colon M \to G$ be the corresponding embedding. Now, let $\sigma=\Sigma \restriction_{im M}$ and $\tau=T \restriction_{im M}$, the compositions of $\Sigma$ and $T$ with the inverse of $\varepsilon$ (which makes sense since the image of $\Sigma$ and $T$ is a manifold in which $\varepsilon$ has an inverse defined). Let $\imath=\Upsilon$. Notice that $\sigma$ and $\tau$, being restrictions to the images of smooth constant rank maps, and the image being an embedded manifold, are automatically surjective submersions. Item $i)$ is just a restatement of item $ii)$ in definition \ref{groupoid_dfn}, restricted to the image. Item $ii)$ is trivial from the definitions. Item $iii)$ follows from
\begin{displaymath}
\sigma(\varepsilon(x))=\Sigma(T(g))=T(g)=\tau(\varepsilon(x))=x
\end{displaymath}
since $g=\varepsilon(x)$ implies $g \in im(T)$, with $g=T(g)$ and $\Sigma T=T$. Item $iv)$ follows from
\begin{displaymath}
\mu(g,\varepsilon(\sigma(g)))=\mu(g,\Sigma(g))=g
\end{displaymath}
and from
\begin{displaymath}
\mu(\varepsilon(\tau(g)),g)=\mu(T(g),g)=g
\end{displaymath}
Item $v)$ is just a restatement of item $iv)$ of definition \ref{groupoid_dfn} and the conditions of item $i)$ in the same definition.
\end{proof}
\end{thm}

The above result is the reason why I prefer to call Lie groupoids those groupoids in which both total space and base space are smooth manifolds. However, differentiable groupoids, in which the total space has a weaker notion of differential structure, are very important in many branches of mathematics and to get such groupoids with this framework is a matter of choosing the right $C^{\infty}$ category.

\begin{expl}
Let $M$ be a smooth manifold. Let $\Aut_{loc}^{\infty}(M)$ be the space of all local diffeomorphisms between open sets of $M$. Define $\Sigma,T,\Upsilon \colon \Aut_{loc}^{\infty}(M) \to \Aut_{loc}^{\infty}(M)$ given by, for all $f \in \Aut_{loc}^{\infty}(M)$,
\begin{eqnarray*}
\Sigma(f) & = & f^{-1} \circ f \\
T(f) & = & f \circ f^{-1} \\
\Upsilon(f) & = & f^{-1}
\end{eqnarray*}
Whenever $\Sigma(g)=T(f)$, define $\mu(g,f)=g \circ f$. We claim that \linebreak$(\Aut_{loc}^{\infty}(M),\Sigma,T,\mu,\Upsilon)$ is a diffeological groupoid, \emph{i.e.} a groupoid in the category of diffeological maps between diffeological spaces. That the space $\Aut_{loc}^{\infty}(M)$ is diffeological, as well as the fact that the maps above defined are diffeological, are left to the reader. The reader unfamiliar with the concept of diffeological space can find a definition in \cite{schreiber_parallel}, but it will not be crucial here, since it only serves as a ``receptacle'' that organizes entities which have some notion of differentiability, such as foliations, singular spaces and spaces of (smooth) functions, just to draw a few examples. The important thing is that this is a category with wild objects, but with very nice categorial properties, allowing the construction of all universal objects needed here.

For the structure, notice that, for all $f \in \Aut_{loc}^{\infty}(M)$,
\begin{displaymath}
T(\Sigma(f))=T(f^{-1} \circ f)=T(id_{dom(f)})=id_{dom(f)}=\Sigma(f)
\end{displaymath}
hence, $T \Sigma=\Sigma$. An analogous reasoning shows that $\Sigma T=T$. It is now easy to see that $\Sigma(\mu(g,f))=\Sigma(f)$ and $T(\mu(g,f))=T(g)$ as well as the universality of $\eta_{\Sigma}$ and $\eta_{T}$ given by $\eta_{\Sigma}(f)=(f,\Sigma(f))$ and $\eta_{T}(f)=(T(f),f)$, for all $f \in \Aut_{loc}^{\infty}(M)$. Besides, $\mu \circ \eta_{\Sigma}=id_{\Aut_{loc}^{\infty}(M)}$ and $\mu \circ \eta_{T}=id_{\Aut_{loc}^{\infty}(M)}$ are also very clear. The associativity of $\mu$ follows directly from the associativity of the composition of functions. For the coequalizer of $(\Sigma,id_{\Aut_{loc}^{\infty}(M)})$, define $\mathcal{O}(M)$ as the space of all open sets of $M$ and $\sigma \colon \Aut_{loc}^{\infty}(M) \to \mathcal{O}(M)$ as the map given by, for all $f \in \Aut_{loc}^{\infty}(M)$,
\begin{displaymath}
\sigma(f)=dom(f)
\end{displaymath}
That $\sigma$ satisfies the desired universal properties is an easy calculation and is left to the reader. The definition of the coequalizer $\tau \colon \Aut_{loc}^{\infty}(M) \to \mathcal{O}(M)$ of the pair $(T,id_{\Aut_{loc}^{\infty}(M)})$ is analogous and left to the reader as well.

One could have started with the space of homeomorphisms between open sets of the topological space $M$, which will be denoted by $\Aut_{loc}^{0}(M)$, getting a similar structure of (topological) groupoid on $\Aut_{loc}^{0}(M)$. In the case when $M$ is a smooth manifold $\Aut_{loc}^{\infty}(M)$ is a wide subgroupoid of $\Aut_{loc}^{0}(M)$, \emph{i.e.} a subgroupoid which contains all the identities of the ambient groupoid. Wide subgroupoids of $\Aut_{loc}^{0}(M)$ play a central role in the theory of manifolds, particularly on the theory of symmetries on manifolds with structure. In fact, a wide subgroupoid $\Gamma \subseteq \Aut_{loc}^{0}(M)$ is called a \emph{pseudogroup of transformations} if it is closed by restrictions, \emph{i.e.} $(f \colon U \to V) \in \Gamma$ if and only if $f_{\alpha}=f|_{U_{\alpha}} \in \Gamma$ for all open covers $U=\bigcup U_{\alpha}$. This definition is equivalent to the one usually found in the literature, for example as in \cite{kobayashi_nomizu}. $\Diamond$
\end{expl}
\begin{rmk}
In the literature on Lie groupoids is often stressed out that a Lie groupoid is a category in which both space of objects and space of morphisms have smooth structures, where all structure maps are smooth, and the source and target (domain and codomain) maps are surjective submersions. The later demand seems as crucial as \emph{ad hoc} and, within this context, it is not clear what must be the conditions in such maps if one attempts to generalise the concept of groupoid for a weaker notion of differentiability, as diffeological spaces. The approach followed in these notes makes clear some things that are more or less obscured by the traditional (in the sense of wide spread) definition. First of all, the category dependent conditions usually imposed on source and target maps (submersion is a feature of smooth spaces) are demonstrable here. Thus, in the category of diffeological spaces of the last example, we didn't need to care which conditions are to be imposed on such maps, avoiding what could be a hard analytical problem. Secondly, this approach makes clear that the ``section of unity'', also stressed out as necessary, is in fact redundant. $\Diamond$
\end{rmk}
\begin{expl}
Another interesting example of these ideas is the groupoid object in $\categories{RAss}$, the category of associative $R$-algebras with its homomorphisms as morphisms. Let $(G,\Sigma,T,\Upsilon,\mu)$ be a groupoid in associative $R$-algebras. The defining conditions of $\Sigma$ and $T$ imply that both are projections. Thus, we have two splitting of $G$
\begin{equation*}
G=ker(\Sigma)\oplus im(\Sigma)
\end{equation*}
\begin{equation*}
G=ker(T)\oplus im(T)
\end{equation*}
The reader can check this by noting that there is a unique way to write an element $g \in G$ in terms of $ker(P)$ and $im(P)$, $P$ being either $\Sigma$ or $T$ (or any projection): $g=(g-P(g))+P(g)$.

The compatibility condition between $\Sigma$ and $T$ implies $im(\Sigma)=im(T)$. It may not be difficult for the reader to check that $H=im(\Sigma)=im(T)$ is a subalgebra of $G$, such that both $ker(\Sigma)$ and $ker(T)$ are $H$-bimodules. Besides, any $x \in H$ satisfies $x=\Sigma(x)=T(x)$. The defining conditions of $\Upsilon$ imply that it is an involution which leaves $H$ invariant. Also, $\Upsilon$ ``reflects'' $ker(\Sigma)$ and $ker(T)$ one to the other.

The product $G \times G$ is the direct product of $R$-algebras, with multiplication defined by
\begin{displaymath}
(g,f)\cdot (g',f')=(gg',ff'),\quad \forall (g,f),(g',f') \in G \times G
\end{displaymath}
Hence, $G_{2}$ is defined as being
\begin{displaymath}
G_{2}=\{ (g,f)\in G \times G \;|\; T(f)=\Sigma(g) \}
\end{displaymath}
Write $G \times G=(ker(\Sigma)\oplus H) \times (ker(T)\oplus H)$. We have then
\begin{eqnarray*}
T(f) = \Sigma(g) & \Rightarrow & T(f_{T}+f_{H})=\Sigma(g_{\Sigma}+g_{H}) \Rightarrow T(f_{H})=\Sigma(g_{H}) \Rightarrow \\
& \Rightarrow & f_{H}=g_{H}
\end{eqnarray*}
where the sub-indices denote the respective decomposition in the direct sum. Hence, $G_{2} \approx ker(\Sigma)\oplus ker(T)\oplus H$. For the record, the induced multiplication in $G_{2}$ reads
\begin{equation*}
(g_{\Sigma},f_{T},x)\cdot(g'_{\Sigma},f'_{T},x')=(g_{\Sigma}g'_{\Sigma}+g_{\Sigma}x'+xg'_{\Sigma},f_{T}f'_{T}+f_{T}x'+xf'_{T},xx')
\end{equation*}
From equation \ref{eta_S_mu_def} in definition \ref{category_dfn}, we have
\begin{displaymath}
\mu(\eta_{\Sigma}(g_{\Sigma}+x))=g_{\Sigma}+x
\end{displaymath}
Which gives
\begin{equation}
\label{mu_decomp_sigma_eqn}
g_{\Sigma}+x=\mu(g_{\Sigma},0,x)
\end{equation}
And for $T$
\begin{displaymath}
\mu(\eta_{T}(g_{T}+x))=g_{T}+x
\end{displaymath}
Which gives
\begin{equation}
\label{mu_decomp_tau_eqn}
g_{T}+x=\mu(0,g_{T},x)
\end{equation}
From these equations we infer that
\begin{enumerate}
\item[i)] $\mu(0,0,x)=x$;
\item[ii)] $\mu(g_{\Sigma},0,0)=g_{\Sigma}$;
\item[iii)] $\mu(0,g_{T},0)=g_{T}$.
\end{enumerate}
But $\mu$, being $R$-linear, must satisfy
\begin{displaymath}
\mu(g_{\Sigma},f_{T},x)=\mu(g_{\Sigma},0,0)+\mu(0,f_{T},0)+\mu(0,0,x)=g_{\Sigma}+f_{T}+x
\end{displaymath}
Hence, $\mu$ is nothing but the sum in the algebra.

From equations \ref{mu_upsilon_eqn} in definition \ref{groupoid_dfn} we have
\begin{eqnarray}
& & \mu(id_{G}\times \Upsilon)(g_{T}+x)=T(g_{T}+x)=x \therefore \nonumber \\
& & x=\mu(g_{T},\Upsilon(g_{T}),x)=g_{T}+\Upsilon(g_{T})+x \therefore \nonumber \\
& & \Upsilon(g_{T})=-g_{T} \label{upsilon1}
\end{eqnarray}
\begin{eqnarray}
& & \mu(id_{G}\times \Upsilon)(g_{\Sigma}+x)=T(g_{\Sigma})+x \therefore \nonumber \\
& & T(g_{\Sigma})+x=\mu(g_{\Sigma},\Upsilon(g_{\Sigma}),x)=g_{\Sigma}+\Upsilon(g_{\Sigma})+x \therefore \nonumber \\
& & \Upsilon(g_{\Sigma})=T(g_{\Sigma})-g_{\Sigma} \label{upsilon2}
\end{eqnarray}
\begin{eqnarray}
& & \mu(\Upsilon \times id_{G})(g_{T}+x)=\Sigma(g_{T})+x \therefore \nonumber \\
& & \Sigma(g_{T})+x=\mu(\Upsilon(g_{T}),g_{T},x)=\Upsilon(g_{T})+g_{T}+x \therefore \nonumber \\
& & \Upsilon(g_{T})=\Sigma(g_{T})-g_{T} \label{upsilon3}
\end{eqnarray}
\begin{eqnarray}
& & \mu(\Upsilon \times id_{G})(g_{\Sigma}+x)=x \therefore \nonumber \\
& & x=\mu(\Upsilon(g_{\Sigma}),g_{\Sigma},x)=\Upsilon(g_{\Sigma})+g_{\Sigma}+x \therefore \nonumber \\
& & \Upsilon(g_{\Sigma})=-g_{\Sigma} \label{upsilon4}
\end{eqnarray}
From \ref{upsilon1} and \ref{upsilon4}, we have that $\Upsilon$ is the additive inversion on the $H$-modules $ker(\Sigma)$ and $ker(T)$. Combining these with \ref{upsilon2} and \ref{upsilon3}, we have
\begin{eqnarray*}
\Upsilon(g_{\Sigma}) & = & T(g_{\Sigma})-g_{\Sigma} \;\therefore \\
-g_{\Sigma} & = & T(g_{\Sigma})-g_{\Sigma} \;\therefore \\
T(g_{\Sigma}) & = & 0
\end{eqnarray*}
An analogous reasoning shows that
\begin{displaymath}
\Sigma(g_{T})=0
\end{displaymath}
Therefore, for all $g \in G$, we have
\begin{eqnarray*}
T(g)-\Sigma(g) & = & T(g_{\Sigma}+x)-\Sigma(g_{\Sigma}+x)=T(g_{\Sigma})+x-\Sigma(g_{\Sigma})-x = \\
& = & 0
\end{eqnarray*}
Thus, $T=\Sigma$.

Let $N=ker(\Sigma)=ker(T)$. We have $G_{2}=N\oplus N\oplus H$. For $p,q \in N$, the multiplication in $G$ must satisfy
\begin{eqnarray*}
p\cdot q & = & \mu(p,0,0)\cdot \mu(0,q,0)=\mu((p,0,0)\cdot (0,q,0)) = \mu(0.p,0.q,0) = \\
& = & 0
\end{eqnarray*}
where we used the fact that $\mu$ is an algebra homomorphism.

It follows that $G$ is the Abelian extension\footnote{For a definition, see for example \cite{loday_cyclic_homology}.} of the algebra $H$ by the $H$-bimodule $N$. The reader may enjoy to prove that any such extension defines a groupoid object in $\categories{RAss}$. The space $\Iso(G,G')$ of isomorphisms of the groupoids $G$ and $G'$ consists of all isomorphic Abelian extensions $G$ and $G'$ of isomorphic $R$-algebras $H$ and $H'$ by the respective (isomorphic) bimodules $N$ and $N'$. In this way one can consider the space $\Hom(G,G')$ of homomorphisms of groupoid objects. It will list all the possible relations (by morphisms) between Abelian extensions. The industrious reader may enjoy to prove in addition that groupoid objects in the categories of Lie algebras and Poisson algebras give Abelian extensions as well. The proof is almost the same of that given here for associative algebras. Hence, we can regard groupoid objects as a kind of trivial extension of the space in question. In cases in which such extensions furnish a representation (as in the algebraic cases above) the groupoid objects are related to classification problems and can offer an organized way to state Morita equivalences. $\Diamond$
\end{expl}

To end this section, let's return to the general situation and see how the definition of groupoid in a category, presented here, has the definition of group in a category as a special case.

\begin{prop}
\label{group_object_prop}
Let $\categories{C}$ be a category with terminal object $1$ and $(G,\Sigma,T,\Upsilon,\mu)$ be a groupoid in $\categories{C}$. Then the groupoid mentioned is a group object in $\categories{C}$, precisely when there is an element $e \colon 1 \to G$ of $G$ such that $\Sigma=T=e_{1}$, where $e_{1}=e\; \circ\; !$, and $!$ denotes the unique morphism $! \colon G \to 1$.
\begin{proof}
Suppose that $\Sigma=T=e\; \circ\; !$. Then, the pull-back $\Pb{G}{\Sigma}{T}{G}$ is the equalizer of the two equal morphisms $\Sigma$ and $T$, which can be written as the composite of the arrow $e$ with the pull-back diagram $\Pb{G}{!}{!}{G}$, furnishing $\Pb{G}{\Sigma}{T}{G}=G \times G$. Thus, $\mu \colon G \times G \to G$. Now we can use the factorization mentioned in remark \ref{product_delta_rmk}, which gives the following diagram for identities:
\begin{displaymath}
\xymatrix@C+20pt{G \times G \ar[r]^-{(id_{G},\Sigma)} & G \times G \ar[d]^{\mu} \\
G \ar[u]_{\Delta} \ar[r]_{id_{G}} & G}
\end{displaymath}
The reader can check the following. There is a canonical isomorphism $\lambda \colon G \to G \times 1$, and this isomorphism can be written as $\lambda=(id_{G},!) \circ \Delta$ (this is a simple exercise in products and terminal objects). The above diagram is then
\begin{displaymath}
\xymatrix@C+20pt{G \times 1 \ar[r]^-{(id_{G},e)} & G \times G \ar[d]^{\mu} \\
G \ar[u]_{\lambda} \ar[r]_{id_{G}} & G}
\end{displaymath}
It is clear that one can do the same for $1 \times G$. For the inversion, we have
\begin{displaymath}
\xymatrix@C+20pt{G \times G \ar[r]^-{(id_{G},\Upsilon)} & G \times G \ar[d]^{\mu} \\
G \ar[u]_{\Delta} \ar[r]_{T} & G}
\end{displaymath}
which gives
\begin{displaymath}
\xymatrix{G \times G \ar[rr]^-{(id_{G},\Upsilon)} & & G \times G \ar[d]^{\mu} \\
G \ar[u]_{\Delta} \ar[r]_{!} & 1 \ar[r]_{e} & G}
\end{displaymath}
It is remaining only the associativity of $\mu$. But this is trivial, because the limit $G_{3}$ is clearly isomorphic to $G \times (G \times G)$ and to $(G \times G) \times G$. These furnish all diagrams in the definition of a group in $\categories{C}$ (see \cite{maclane_categories}).

On the other hand, if one has a group object in $\categories{C}$, a groupoid is obtained by defining $\Sigma=T=e\;\circ\;!$.
\end{proof}
\end{prop}

\section{Cogroupoids}
\label{cogroupoids_sec}
\begin{dfn}[Cogroupoid in a category]
\label{cogroupoid_dfn}
Let $\categories{C}$ be a category. A cogroupoid in $\categories{C}$ is a groupooid in $\categories{C}^{op}$, where $\categories{C}^{op}$ denotes the opposite category of $\categories{C}$.
\end{dfn}
It is worthy to write explicitly the diagrams in the case in which the given colimits can be expressed as combinations of coproducts and coequalizers. In this case a cogroupoid in $\categories{C}$ consists in $(\mathcal{C},\mathcal{S},\mathcal{T},\mathcal{U},m)$, where $\mathcal{C} \in \categories{C}_{0}$, $\mathcal{S},\mathcal{T},\mathcal{U} \in \End_{\categories{C}}(\mathcal{C})$ and $m \in \categories{C}(\mathcal{C},\mathcal{C}^{2})$, with $\mathcal{C}^{2} \in \categories{C}_{0}$ being the push-out
\begin{displaymath}
\xymatrix{\mathcal{C} \ar[d]_{\mathcal{S}} \ar[r]^{\mathcal{T}} & \mathcal{C} \ar[d]^{\imath_{2}} \\
\mathcal{C} \ar[r]_{\imath_{1}} & \mathcal{C}^{2}}
\end{displaymath}
We also have $\mathcal{S}^{2}=\mathcal{S}$, $\mathcal{T}^{2}=\mathcal{T}$, $\mathcal{S}\mathcal{T}=\mathcal{S}$, $\mathcal{T}\mathcal{S}=\mathcal{T}$, $\mathcal{U}\mathcal{T}=\mathcal{S}$, $\mathcal{S}\mathcal{U}=\mathcal{S}$, $\mathcal{U}^{2}=id_{\mathcal{C}}$. It is commutative
\begin{displaymath}
\xymatrix{\mathcal{C} \ar[d]_{\mathcal{T}} \ar[r]^{\mathcal{T}} & \mathcal{C} \ar[d]^{m} & \mathcal{C} \ar[l]_{\mathcal{S}} \ar[d]^{\mathcal{S}} \\
\mathcal{C} \ar[r]_{\imath_{1}} & \mathcal{C}^{2} & \mathcal{C} \ar[l]^{\imath_{2}}}
\end{displaymath}
$\mathcal{C} \sqcup \mathcal{C},\mathcal{C} \sqcup \mathcal{C} \sqcup \mathcal{C} \in \categories{C}_{0}$ and the codiagonal morphism $\delta \colon \mathcal{C} \sqcup \mathcal{C} \to \mathcal{C}$ do exists, making commutative the following diagram
\begin{displaymath}
\xymatrix@C+20pt{\mathcal{C} \ar[d]_{m} \ar[r]^{id_{\mathcal{C}}} & \mathcal{C} & \mathcal{C} \ar[l]_{id_{\mathcal{C}}} \ar[d]^{m} \\
\mathcal{C}^{2} \ar[r]_-{\mathcal{T} \sqcup id_{\mathcal{C}}} & \mathcal{C} \sqcup \mathcal{C} \ar[u]_{\delta} & \mathcal{C}^{2} \ar[l]^-{id_{\mathcal{C}} \sqcup \mathcal{S}}}
\end{displaymath}
where we interpret symbols like $\mathcal{T} \sqcup id_{\mathcal{C}}$ as the morphism induced by the epic $\mathcal{C} \sqcup \mathcal{C} \to \mathcal{C}^{2} $, since push-outs are ``quotient objects'' of coproducts (by coequalizers). It is commutative
\begin{displaymath}
\xymatrix@C+20pt{\mathcal{C} \ar[d]_{m} \ar[r]^{\mathcal{S}} & \mathcal{C} & \mathcal{C} \ar[l]_{\mathcal{T}} \ar[d]^{m} \\
\mathcal{C}^{2} \ar[r]_-{\mathcal{U} \sqcup id_{\mathcal{C}}} & \mathcal{C} \sqcup \mathcal{C} \ar[u]_{\delta} & \mathcal{C}^{2} \ar[l]^-{id_{\mathcal{C}} \sqcup \mathcal{U}}}
\end{displaymath}
$\mathcal{C}^{3} \in \categories{C}_{0}$ is the evident colimit that makes the following diagram commutative
\begin{displaymath}
\xymatrix@C+20pt{\mathcal{C} \ar[d]_{m} \ar[r]^{m} & \mathcal{C}^{2} \ar[d]^{id_{\mathcal{C}} \sqcup m} \\
\mathcal{C}^{2} \ar[r]_-{m \sqcup id_{\mathcal{C}}} & \mathcal{C}^{3}}
\end{displaymath}
and the pairs $(\mathcal{S},id_{\mathcal{C}})$ and $(\mathcal{T},id_{\mathcal{C}})$ have equalizers.
\begin{rmk}
When the coproduct $\mathcal{C}\sqcup \mathcal{C}$ is defined, the codiagonal morphism is canonically defined as the unique morphism $\delta \colon \mathcal{C}\sqcup \mathcal{C} \to \mathcal{C}$ making commutative the following diagram
\begin{displaymath}
\xymatrix{ & \mathcal{C} & \\
& \mathcal{C} \sqcup \mathcal{C} \ar[u]^{\delta} & \\
\mathcal{C} \ar[ur]_{\imath_{1}} \ar@/^/[uur]^{id_{\mathcal{C}}} & & \mathcal{C} \ar[ul]^{\imath_{2}} \ar@/_/[uul]_{id_{\mathcal{C}}}}
\end{displaymath}
As an example, if $A$ is an associative commutative unital $R$-algebra, the coproduct is the tensor product\footnote{See, for example, \cite{lang_algebra}.} (over $R$) and the diagram means $\delta(a\otimes 1)=a$ and $\delta(1\otimes a)=a$, for all $a \in A$. But the multiplication in $A$ has this property. Being $\delta$ unique, $\delta$ is just the multiplication in $A$. $\Diamond$
\end{rmk}
\begin{expl}
In \cite{nestruev_smooth_manifolds}, the authors show that there is a category $C^{\infty}\categories{Alg}$, full subcategory of the category of associative commutative unital $\R$-algebras, which is (anti)isomorphic to $\categories{Man}$. The functor which gives the isomorphism from manifolds to $C^{\infty}\categories{Alg}$ consists in taking the $\R$-algebra of smooth functions on a manifold (on objects) and in taking the homomorphism of $\R$-algebras induced by composition (on morphisms). Thus, our analysis shows that cogroupoids in $C^{\infty}\categories{Alg}$ are exactly the Lie groupoids in $\categories{Man}$. $\Diamond$
\end{expl}
In light of proposition \ref{group_object_prop} it is worthy to give a description when the cogroupoid is dual of a groupoid that is a group object, in a category which is subcategory of associative commutative unital $R$-algebras (for $R$ a ring).
\begin{prop}
\label{hopf_algebra_object_prop}
Let $R$ be a commutative unital ring, $\categories{C}$ be the category of associative commutative unital $R$-algebras (and homomorphisms) and \linebreak$(\mathcal{C},\mathcal{S},\mathcal{T},\mathcal{U},m)$ be a cogroupoid in $\categories{C}$. Then the cogroupoid mentioned is a Hopf algebra, precisely when there is a homomorphism of $R$-algebras $\varepsilon \colon \mathcal{C} \to R$, such that $\mathcal{S}=\mathcal{T}= \imath \circ \varepsilon$, where $\imath \colon R \to \mathcal{C}$ is the canonical inclusion of $R$ in the $R$-algebra $\mathcal{C}$, satisfying $\varepsilon \circ \imath=id_{R}$.
\begin{proof}
Let's begin by constructing $\mathcal{C}^{2}$. The push-out of two $R$-algebras of the given type is obtained by coequalizing the coproduct of $\mathcal{C}$ and $\mathcal{C}$ by the homomorphisms $\mathcal{S}$ and $\mathcal{T}$. In other words, one just have to take the tensor algebra $\mathcal{C} \otimes \mathcal{C}$ (tensor over $R$) and then take the quotient by the ideal generated by elements of the form $\mathcal{S}(c)\otimes 1-1 \otimes \mathcal{T}(c)$, for all $c \in \mathcal{C}$, where $1$ denotes the unit element in $\mathcal{C}$. By the hypothesis, we have $\mathcal{S}(c)=\mathcal{T}(c)$, for all $c \in \mathcal{C}$. Hence $\mathcal{C}^{2}=\mathcal{C} \otimes \mathcal{C}$. It is not hard to see that $\mathcal{C}^{3}$ is isomorphic to both $(\mathcal{C} \otimes \mathcal{C}) \otimes \mathcal{C}$ and $\mathcal{C} \otimes (\mathcal{C} \otimes \mathcal{C})$.  By hypothesis, $m \colon \mathcal{C} \to \mathcal{C} \otimes \mathcal{C}$ is a homomorphism of $R$-algebras and it is associative. By calling $\mu \colon \mathcal{C} \otimes \mathcal{C} \to \mathcal{C}$ the multiplication on $\mathcal{C}$ (and remember that it is also the codiagonal morphism), we have $(\mathcal{C},\mu,\imath)$ an $R$-algebra, $(\mathcal{C},m,\varepsilon)$ an $R$-coalgebra, such that $(\mathcal{C},\mu,m,\imath,\varepsilon)$ is an $R$-bialgebra. Now, the diagram
\begin{displaymath}
\xymatrix{\mathcal{C} \ar[d]_{m} \ar[r]^{\varepsilon} & R \ar[r]^{\imath} & \mathcal{C} & R \ar[l]_{\imath} & \mathcal{C} \ar[l]_{\varepsilon} \ar[d]^{m} \\
\mathcal{C} \otimes \mathcal{C} \ar[rr]_-{\mathcal{U} \otimes id_{\mathcal{C}}} & & \mathcal{C} \otimes \mathcal{C} \ar[u]_{\mu} & & \mathcal{C} \otimes \mathcal{C} \ar[ll]^-{id_{\mathcal{C}} \otimes \mathcal{U}}}
\end{displaymath}
is just the diagram of an antipode on $\mathcal{C}$, making $(\mathcal{C},\mu,m,\imath,\varepsilon,\mathcal{U})$ a Hopf algebra. As a special feature, from the diagrams involving $m$ and $\mathcal{S}$ one can easily see that $m(1)=1 \otimes 1$.

On the other hand, given a Hopf $R$-algebra $(\mathcal{C},\mu,m,\imath,\varepsilon,\mathcal{U})$, define $\mathcal{S}=\mathcal{T}=\imath \circ \varepsilon$. It is easy to see that $(\mathcal{C},\mathcal{S},\mathcal{T},\mathcal{U},m)$ is a cogroupoid in $\categories{C}$ of the given type.
\end{proof}
\end{prop}
The proposition above legitimate the terminology \emph{Hopf algebroid} for a cogroupoid object in a category of algebras (with some ``nice'' properties), largely used in the literature.

\section{Category actions}
\label{groupoid_actions_sec}
It may be interesting at this point to give a diagrammatic description of an action of a category object in a category on an object. It will be described left actions, but one can easily define right actions as well. The only thing to be noted is that the situation here is not as symmetric as in the case of a group, for example. For a right action, the role of $\Sigma$ and $T$ are interchanged.  Here $\categories{C}$ denotes a category.
\begin{dfn}[Category action]
\label{groupoid_act_dfn}
Let $(C,\Sigma,T,\mu)$ be a category in $\categories{C}$, $E \in \categories{C}_{0}$ be an object and $\varphi \in \categories{C}(E,C)$ be a morphism. Suppose that the pull-back $\Pb{C}{\Sigma}{\varphi}{E}$ given by
\begin{displaymath}
\xymatrix{ & & \Pb{C}{\Sigma}{\varphi}{E} \ar[dl]_-{\pi_{\Sigma \varphi 1}} \ar[dr]^-{\pi_{\Sigma \varphi 2}} & & \\
& C \ar[dl]_-{T} \ar[dr]_-{\Sigma} & & E \ar[dl]^-{\varphi_{\Sigma}} & \\
C & & C & & }
\end{displaymath}
where $\varphi_{\Sigma}=\Sigma \circ \varphi$, do exists. A morphism $\theta \colon \Pb{C}{\Sigma}{\varphi}{E} \to E$ in $\categories{C}$ is called an action of $C$ on $E$ by $\varphi$ if and only if it satisfies:
\begin{enumerate}
\item[i)] by denoting $\eta_{\Sigma \varphi} \colon C \to \Pb{C}{\Sigma}{\varphi}{E}$ as the unique arrow given by the cone $C \xleftarrow{\varphi_{\Sigma}} E \xrightarrow{id_{E}} E$, such that $\pi_{\Sigma \varphi 1} \circ \eta_{\Sigma \varphi}=\varphi_{\Sigma}$ and $\pi_{\Sigma \varphi 2} \circ \eta_{\Sigma \varphi}=id_{E}$, then it commutes
\begin{equation}
\label{act_tau_eqn}
\xymatrix@C+40pt{ & E \ar[dl]_-{\varphi_{\Sigma}} \ar[d]^-{\eta_{\Sigma \varphi}} \ar@/^3pc/[dd]^-{id_{E}} \\
C & \Pb{C}{\Sigma}{\varphi}{E} \ar[l]_-{T \circ \pi_{\Sigma \varphi 1}} \ar[d]^-{\theta} \\
& E \ar[ul]^-{\varphi_{\Sigma}}}
\end{equation}
\item[ii)] there exist the pull-back $\Pb{(\Pb{C}{\Sigma}{T}{C})}{\Sigma_{2}}{\varphi}{E}$ given by
\begin{equation}
\label{G2xE_def}
\xymatrix{
& & \Pb{(\Pb{C}{\Sigma}{T}{C})}{\Sigma_{2}}{\varphi}{E} \ar[dl]_-{\pi_{\Sigma_{2}\varphi 1}} \ar[dr]^-{\pi_{\Sigma_{2}\varphi 2}} & & \\
& \Pb{C}{\Sigma}{T}{C} \ar[dl]_-{T \circ \pi_{\Sigma T1}} \ar[dr]^-{\Sigma \circ \pi_{\Sigma T2}} & & E \ar[dl]_-{\varphi_{\Sigma}} & \\
C & & C & & }
\end{equation}
and the pull-back $\Pb{C}{\Sigma}{T_{2}}{(\Pb{C}{\Sigma}{\varphi}{E})}$ given by a similar diagram;
\item[iii)] by denoting $\mu \times id_{E} \colon \Pb{(\Pb{C}{\Sigma}{T}{C})}{\Sigma_{2}}{\varphi}{E} \to \Pb{C}{\Sigma}{\varphi}{E}$ as the unique arrow given by the cone $C \xleftarrow{\mu \circ \pi_{\Sigma_{2}\varphi 1}} \Pb{(\Pb{C}{\Sigma}{T}{C})}{\Sigma_{2}}{\varphi}{E} \xrightarrow{id_{E} \circ \pi_{\Sigma_{2}\varphi 2}} E$ and by $id_{C} \times \theta \colon \Pb{C}{\Sigma}{T_{2}}{(\Pb{C}{\Sigma}{\varphi}{E})} \to \Pb{C}{\Sigma}{\varphi}{E}$ as the unique arrow given by the cone $C \xleftarrow{id_{C} \circ \pi_{\Sigma T_{2},1}} \Pb{C}{\Sigma}{T_{2}}{(\Pb{C}{\Sigma}{\varphi}{E})} \xrightarrow{\theta \circ \pi_{\Sigma T_{2},2}} E$, it holds true
\begin{equation}
\label{act_mu_def}
\theta \circ (id_{C} \times \theta) \approx \theta \circ (\mu \times id_{E})
\end{equation}
meaning that the left hand side of the last equation is the same as the right hand side up to isomorphism. Again, this is necessary because we have in general only an isomorphism between $\Pb{(\Pb{C}{\Sigma}{T}{C})}{\Sigma_{2}}{\varphi}{E}$ and $\Pb{C}{\Sigma}{T_{2}}{(\Pb{C}{\Sigma}{\varphi}{E})}$.
\end{enumerate}
\end{dfn}

One important fact is that a category internal to a category acts on itself in a natural way. To see this, take $\varphi=T$ and $\theta=\mu$. Now, notice that $T_{\Sigma}=\Sigma \circ T=T$, furnishing $\eta_{\Sigma T}=\eta_{T}$ and fulfilling the upper triangle in diagram \ref{act_tau_eqn}. The remaining conditions in diagram \ref{act_tau_eqn} are fulfilled, because of conditions \ref{mu_tau_sigma_def} and \ref{eta_S_mu_def} in definition \ref{category_dfn}. Diagram \ref{act_mu_def} is just a statement of the associativity of $\mu$.

With these tools at hand, consider again the case of groupoids in $\categories{Man}$, \emph{i.e.} Lie groupoids. Associated to the projection $T$ there is the vertical space $V^{T}G=ker(\T T)$. The multiplication $\mu$ of the groupoid $G$ induces a left action $L_{g} \colon G_{\Sigma(g)} \to \; {}_{T(g)} \! G$, where $G_{\Sigma(g)}$ denotes the subset of elements $f \in G$ such that $L_{g}(f)=\mu(g,f)$ makes sense, and ${}_{T(g)} \! G=\{ h \in G \;|\; T(h)=T(g) \}$. Let $\mathfrak{X}_{LI}^{T}(G)$ be the subspace of $\Gamma(V^{T}G)$ given by the left invariant vector fields, \emph{i.e.} vector fields which satisfy
\begin{displaymath}
\T_{f}L_{g}(X_{f})=X_{\mu(g,f)} \quad \forall (g,f) \in G_{2}.
\end{displaymath}
By denoting $AG$ the restriction of the tangent bundle $\T G$ to points of $M=im(T)=im(\Sigma)$, it is easy to see that there is an isomorphism between $\mathfrak{X}_{LI}^{T}(G)$ and $\Gamma(AG)$ given by $\Phi \colon \mathfrak{X}_{LI}^{T}(G) \to \Gamma(AG)$, $\Phi(X)=X \circ \Sigma$, since left invariant vector fields are determined by its values at $M$ (see \cite{mackenzie_lie_groupoids}, for details). Let $C^{\infty}_{\Sigma I}(G)$ be the $\Sigma$-invariant smooth functions, \emph{i.e.} the set of $\eta \in C^{\infty}(G)$ such that $\eta=\eta \circ \Sigma$. It is also easy to see that $C^{\infty}_{\Sigma I}(G)$ is isomorphic to $C^{\infty}(M)$. Let $a=\Gamma(\T \Sigma)$. For all $X \in \mathfrak{X}_{LI}^{T}(G)$ (regarded as an element of $\Gamma(AG)$) and for all $\eta \in C^{\infty}_{\Sigma I}(G)$, we have
\begin{displaymath}
\T_{g} \Sigma(X_{g})(\eta)=X_{\Sigma(g)}(\eta \circ \Sigma)=X_{\Sigma(g)}(\eta)
\end{displaymath}
for all $g \in G$. $\Gamma(AG)$ is closed for the Lie bracket, hence $\mathfrak{X}_{LI}^{T}(G)$ is also closed. This fact with the later equation gives
\begin{equation}
\label{anchor_Lie_hom_eqn}
a(\liebr{X}{Y})=\liebr{a(X)}{a(Y)} \quad \forall X,Y \in \mathfrak{X}_{LI}^{T}(G)
\end{equation}
It is clear that $\mathfrak{X}_{LI}^{T}(G)$ is a $C^{\infty}_{\Sigma I}(G)$-module. For all $X,Y \in \mathfrak{X}_{LI}^{T}(G)$ and for all $\eta, \theta \in C^{\infty}_{\Sigma I}(G)$, we have
\begin{eqnarray*}
\liebr{X}{\eta Y}(\theta) & = & X(\eta Y(\theta))-\eta Y(X(\theta)) = \\
& = & X(\eta)Y(\theta)+\eta X(Y(\theta))-\eta Y(X(\theta)) = \\
& = & \eta(X(Y(\theta))-Y(X(\theta)))+X(\eta \circ \Sigma)Y(\theta) = \\
& = & \eta \liebr{X}{Y}(\theta)+a(X)Y(\theta)
\end{eqnarray*}
which means
\begin{equation}
\label{anchor_Leibniz_eqn}
\liebr{X}{\eta Y}=\eta \liebr{X}{Y}+a(X)Y
\end{equation}
Equations \ref{anchor_Lie_hom_eqn} and \ref{anchor_Leibniz_eqn} plus the module structures and identifications above mentioned make $(\mathfrak{X}_{LI}^{T}(G),C^{\infty}_{\Sigma I}(G))$ into a \emph{Lie-Rinehart pair}. For a definition the reader can see \cite{rinehart_differential_forms}. The issue here is that this is an algebraic construction (or version) of the \emph{Lie algebroid} of the Lie groupoid $G$ (see \cite{mackenzie_lie_groupoids} - there, it was used right actions). The language used here made clear that the Lie algebroid structure is obtained by projecting with $\T \Sigma$, vector fields which are $T$-invariant in some sense.

\section{Conclusions and further directions}
\label{epilogue_sec}
The aim of these notes is to call attention to a technique apparently overlooked. This is not for the substitution of the well established procedures, but for aggregation of a new look that may be useful in some situations (but certainty not in all situations). However, the tools as presented here can be carried out to cover (I hope) all the main concepts involving groupoids. Let me say a few words about the framework.

If $(G,\Sigma,T,\Upsilon,\mu)$ is a groupoid in $\categories{C}$, one can consider the subobject $M$ of $G$ given by the coequalizer conditions. A bisection of $G$ can be defined as a subobject $N \rightarrowtail G$ such that both $\Sigma$ and $T$, when pulled-back to $N$, give $\Sigma^{\ast} \colon M \to N$ and $T^{\ast} \colon M \to N$ as being isomorphisms. This is just an example of how to obtain more complicated structures, often useful in the standard theory.

Dually, one can obtain similar concepts for cogroupoids, as coactions and biretractions. It is the belief of the author that the duality between smooth manifolds and its smooth algebras allows a clean description of jet groupoids by means of algebraic jets of cogroupoids in smooth algebras (in the sense of \cite{nestruev_smooth_manifolds}), since there is a natural duality between jets of algebras and high order differential operators. Thus, one could discuss jet prolongations in Hopf algebroids (at least in the commutative case) and coactions of such jets. This is a study in progress by this author.

Another direction of study is about groupoid objects in Poisson algebras. Such groupoids can describe extensions of Poisson algebras by ideals that could be useful on the understanding of the classical BRST formalism, in particular, how the projections could describe constrains in the classical setting.

As a final remark, the ideas of this work can be useful in the process of enrichment and internalization of categories. Notably this can be useful on the description of $n$-categories and $n$-groupoids. Although straightforward in the strict case, the weak (and more important) case may demand a bit more work. However, I think it is worthy.

\section*{Acknowledgements}
I wish to thanks Dr. Michael Forger, my advisor, for the discussions and for call my attention to this problem, and also, the group of research in Hamiltonian covariant formalism of Instituto de Matematica e Estat\'istica of Universidade de S\~ao Paulo. A special thanks to Sebasti\'an Vidal, for the enlightened discussions. This work was supported by CAPES.

%% file: appendix.tex
\section[Appendix]{Appendix}
\label{appendix_sec}
The first remark is about the category of smooth manifolds. Traditionally this category has as objects those obtained by gluing open subsets of $\R^{n}$, for some $n \geq 1$. This is the idea of the definition of manifolds by charts and atlases. By the embedding theorem of Whitney, this is essentially a small category. However, it lacks a terminal object. This is not a big issue, since it seems natural to consider points of $\R^{n}$ as being $0$-dimensional manifolds. Indeed, one element sets have a unique topology which can be considered with a (trivial) differentiable structure, which in turn furnishes $\R$ as the ring of $C^{\infty}$ functions. This solves the problem of the terminal object and furnishes an important element on the dual: the ground field. Hence the only thing to do is to choose one such set and add it to the category. This is important because if one attempt to add all one point sets, the category obtained is no longer small.

The definition of category used in this paper to give an one object description of a groupoid follows \cite{freyd_categories} and is reproduced here, with some notational modifications, for convenience.
\begin{dfn}[Category]
A category consists of individuals $f,g,h,\ldots$, of unary function symbols $\Sigma,T$ and of a binary partial function symbol $\mu$, satisfying
\begin{enumerate}
\item $\forall f(T \Sigma f = \Sigma f)$;
\item $\forall f(\Sigma T f = T f)$;
\item $\forall g \forall f(\exists h(\mu(g,f) = h) \Leftrightarrow (\Sigma g = T f))$;
\item $\forall f(\mu(f,\Sigma f) = f)$;
\item $\forall f(\mu(T f,f) = f)$;
\item $\forall g\forall f(\exists h(\mu(g,f) = h) \Rightarrow (\Sigma \mu(g,f) = \Sigma f))$;
\item $\forall g\forall f(\exists h(\mu(g,f) = h) \Rightarrow (T \mu(g,f) = T g))$;
\item $\forall h \forall g \forall f((\Sigma h = T g) \wedge (\Sigma g = T f) \Rightarrow (\mu(\mu(h,g),f) = \mu(h,\mu(g,f))))$.
\end{enumerate}
\end{dfn}
The reader who wishes a first order theory may replace the rather awkward ``partial function symbol'' $\mu$ by a ternary predicative letter with the desired properties. The nice side of presenting $\mu$ as a partial function symbol is that it makes clear that it is a quasi-algebraic theory, in the sense that $\mu$ furnishes an equation as long as a condition is fulfilled. In the cited reference the authors call it an \emph{essentially algebraic} theory.

The reader can check that an interpretation of this theory is equivalent to consider a category object in $\categories{Set}$, whatever (reasonable) category of sets is at hand. Here, the coequalizers conditions are missing, simply by the fact that such categories are co-complete (so, the existence of the coequalizers is automatic).

Hence, the one object approach gives the means to build ``interpretations'' of a category theory in a non set theoretic environment, in very general situations, without any additional (\emph{ad hoc}) suppositions.